\documentclass[12pt,leqno]{amsart}
\usepackage{amssymb,amsmath}
\begin{document}

\title[Ramanujan Property for Regular Cubical
Complexes]{The Ramanujan Property for Regular\\
 Cubical Complexes}

\author{Bruce W. Jordan}
\address{Department of Mathematics, Box G-0930\\
Baruch College, CUNY\\
17 Lexington Avenue\\
New York, NY 10010, USA}
\thanks{The first author was partially supported
by grants from the NSF and PSC-CUNY}
\author{Ron Livn\'{e}}
\address{Mathematics Institute\\
The Hebrew University of Jerusalem\\
Givat Ram, Jerusalem  91904  Israel}
\email{rlivne@sunset.ma.huji.ac.il}
\thanks{Both authors were supported by a joint
Binational Israel-USA Foundation grant}
\date{October 4, 1999}
\begin{abstract} We consider cubical complexes
which are uniformized by an ordered product of
regular trees. For these we define the notion
of being {\em Ramanujan}, generalizing the
one-dimensional definition introduced by Lubotzky,
Phillips, and Sarnak \cite{LPS}. As in \cite{JL6},
we also allow local systems. We discuss the
significance of this property, and then we
construct explicit arithmetic examples using
quaternion algebras over totally real fields.
Here we reduce the Ramanujan property to special
cases of the Ramanujan-Petersson conjecture, many
of which are known. Our examples subsume the
constructions of \cite{LPS}, \cite{JL6}, and
\cite{Moz}.
\end{abstract}
\maketitle

\newcounter{soussec}
\newcommand{\soussection}[1]{\par \smallskip
\par \noindent {\bf \arabic{section}.\arabic{soussec}\
#1 \hspace*{1em}} \refstepcounter{soussec}}

\newtheorem{theorem}{Theorem}[section]
\newtheorem{proposition}[theorem]{Proposition}
\newtheorem{variant}[theorem]{A variant}
\newtheorem{fact}[equation]{}
\newtheorem{corollary}[theorem]{Corollary}
\newtheorem{definition}[theorem]{Definition}
\newtheorem{example}[theorem]{Example}
\newtheorem{examples}[theorem]{Examples}
\newtheorem{conditions}[theorem]{Conditions}
\newtheorem{remark}[theorem]{Remark}
\newtheorem{remarks}[theorem]{Remarks}
\newtheorem{lemma}[theorem]{Lemma}
\newtheorem{maintheorem}[theorem]{Main Theorem}

\newenvironment{pf}{\noindent {\it Proof:}\/ }
{\par \par \medskip\par }

\newcommand{\ra}{\rightarrow}
\newcommand{\sra}[1]{\stackrel{#1}{\rightarrow}}
\newcommand{\la}{\leftarrow}
\newcommand{\lra}{\longrightarrow}
\newcommand{\ov}{\overline}

\newcommand{\set}[1]{\{\,#1\,\}}
\newcommand{\Disc}{{\operatorname{Disc}}}
\newcommand{\dist}{{\operatorname{dist}}}
\newcommand{\edg}{{\operatorname{Ed}}}
\newcommand{\edgo}[1]
        {{\operatorname{Ed}^{\rm o}_{#1}}}
\newcommand{\GL}{{\operatorname{GL}}}
\newcommand{\Gal}{{\operatorname{Gal}}}
\newcommand{\Gr}{{\operatorname{Gr}}}
\newcommand{\Hom}{{\operatorname{Hom}}}
\newcommand{\Aut}{{\operatorname{Aut}}}
\newcommand{\Ima}{{\operatorname{Im}}\,}
\newcommand{\Ker}{{\operatorname{Ker}}\,}
\renewcommand{\L}{{\operatorname{L}}}
\newcommand{\Ltwo}{{\L_{2}}}
\newcommand{\Matwo}{{\rm Mat}_{2\times 2}}
\newcommand{\Nm}{{\operatorname{Nm}}}
\newcommand{\PGL}{{\rm PGL}}
\newcommand{\PSL}{{\operatorname{PSL}}}
\newcommand{\res}{\operatorname{res}}
\newcommand{\Res}{{\operatorname{Res}}}
\newcommand{\SL}{{\operatorname{SL}}}
\newcommand{\sign}{{\operatorname{sign}}}
\newcommand{\SU}{{\operatorname{SU}}}
\newcommand{\St}{{\operatorname{St}}}
\newcommand{\Symm}{{\operatorname{Symm}}}
\newcommand{\Tr}{{\operatorname{\Delta}}}
\newcommand{\Tra}{{\operatorname{Tr}}}
\newcommand{\Ver}{{\operatorname{Ver}\,}}
\newcommand{\val}{{\operatorname{val}\,}}
\newcommand{\vk}{{\vec{k}}}
\newcommand{\vr}{{\vec{r}}}
\newcommand{\vs}{{\vec{s}}}
\newcommand{\pard}[1]{\partial_{#1}}
\newcommand{\parad}[1]{\partial^\ast_{#1}}
\newcommand{\der}{{d}}
\newcommand{\derad}{\der^\ast}
\newcommand{\parj}{\pard{j}}
\newcommand{\paradj}{\parad{j}}
\newcommand{\bul}{\bullet}
\newcommand{\Id}{{\rm Id}}
\newcommand{\Iw}{{U}}
\newcommand{\Lap}[1]{\square_{#1}}
\newcommand{\Lapj}{\Lap{j}}
\newcommand{\Lapt}[1]{\Lap{{\rm tot},#1}}
\newcommand{\ocell}[1]{{\Sigma_{#1}}}
\newcommand{\ZtwoI}{(\ZZ/2\ZZ)^I}
\newcommand{\inv}[1]{{\rm inv}_{#1}}
\newcommand{\topf}[1]{{\rm top}_{#1}}
\newcommand{\botf}[1]{{\rm bot}_{#1}}
\newcommand{\invj}{\inv{j}}
\newcommand{\topj}{\topf{j}}
\newcommand{\botj}{\botf{j}}
\newcommand{\topi}{\topf{i}}
\newcommand{\boti}{\botf{i}}
\newcommand{\ktype}{$(k_1,\ldots,k_g)$\/}
\newcommand{\rtype}{$(r_1,\ldots,r_g)$\/}
\newcommand{\HQ}{{\rm H}}
\newcommand{\pib}{{\ov{\pi}}}
\newcommand{\nub}{{\ov{\nu}}}
\newcommand{\gammab}{{\ov{\gamma}}}
\newcommand{\SN}{{\SL_2(\FF_N)}}
\newcommand{\PSN}{{\PSL_2(\FF_N)}}
\newcommand{\phib}{{\ov{\phi}}}
\newcommand{\cMF}{{\cM_{0,F}}}
\newcommand{\leg}[2]{{\left(\frac{#1}{#2}\right)}}

\newcommand{\calA}{{\mathcal A}}
\newcommand{\calB}{{\mathcal B}}
\newcommand{\calH}{{\mathcal H}}
\newcommand{\calI}{{\mathcal I}}
\newcommand{\calJ}{{\mathcal J}}
\newcommand{\calL}{{\mathcal L}}
\newcommand{\calM}{{\mathcal M}}
\newcommand{\calO}{{\mathcal O}}
\newcommand{\calT}{{\mathcal T}}
\newcommand{\calX}{{\mathcal X}}
\newcommand{\cA}{{\calA}}
\newcommand{\cB}{{\calB}}
\newcommand{\cH}{{\calH}}
\newcommand{\cL}{{\calL}}
\newcommand{\cM}{{\calM}}
\newcommand{\cO}{{\calO}}
\newcommand{\cT}{{\calT}}
\newcommand{\cX}{X}
\newcommand{\bG}{{\bf G}}
\newcommand{\bH}{{\bf H}}
\newcommand{\bZ}{{\bf Z}}
\newcommand{\cXtil}{\tilde{\cX}}
\newcommand{\fA}[2]{{#1^{{\rm f}#2}}}
\newcommand{\rf}{{\rm f}}

\renewcommand{\AA}{{\mathbb A}}
\newcommand{\CC}{{\mathbb C}}
\newcommand{\FF}{{\mathbb F}}
\newcommand{\HH}{{\mathbb H}}
\newcommand{\ZZ}{{\mathbb Z}}
\newcommand{\QQ}{{\mathbb Q}}
\newcommand{\RR}{{\mathbb R}}

\newcommand{\sig}{\sigma}
\newcommand{\lam}{\lambda}
\newcommand{\vpi}{\varpi}
\newcommand{\FFq}{\FF_q}

\newcommand{\qn}{{\hat 1}}
\newcommand{\qi}{{\hat\imath}}
\newcommand{\qj}{{\hat\jmath}}
\newcommand{\qk}{{\hat k}}

\newcommand{\ip}[2]{{\langle #1 ,#2\rangle}}
\newcommand{\triv}[1]{{{\mathcal T}_{#1}}}

\newcommand{\matr}[4]
{\left[\begin{array}{cc} #1 & #2 \\ #3 & #4
              \end{array} \right]}
\newcommand{\smatr}[4]{\Bigl[\begin{array}
{@{\hspace{0.1em}}c@{\hspace{0.5em}}c@{\hspace{0.1em}}}
   \scriptstyle #1 & \scriptstyle #2 \\
   \scriptstyle #3 & \scriptstyle #4 \end{array}
   \Bigr]}
\newcommand{\Tptild}{\tilde{T}_p}

\section*{Introduction}
Ramanujan graphs were defined by Lubotzky,
Phillips, and Sarnak in \cite{LPS} as regular
graphs whose adjacency matrices, or their
laplacians, have eigenvalues satisfying some
``best possible'' bounds. Such graphs possess
many interesting properties.  In this paper we
will give a higher dimensional generalization
of this theory to {\em regular cubical
complexes}. By definition, 
$\vr = (r_1,\dots,r_g)$\/-regular complexes
are cell
complexes locally isomorphic to the (ordered)
product of $g$ regular trees, with the $j$\/th
tree of regularity $r_j\geq 3$. Each cell is
an $i$\/-cube (i.e., an $i$\/-dimensional cube)
with $0\leq i\leq g$. Throughout each 
$(g-1)$\/-cube exactly one of the tree factors,
say the $j$\/th one, is constant, and there are
$r_j$ $g$\/-cubes passing through it. When $g=1$ we
simply have an $r$\/-regular graph.

The spaces of $i$\/-cochains $C^i(X)$ (with real
or complex coefficients) of a finite cubical
complex $X$ are inner product vector spaces with
an orthonormal basis corresponding to the
characteristic functions of the $i$\/-cells.
There are partial boundary operators
$\parj=\pard{j,i}:C^i(X)\ra C^{i+1}(X)$ for
$1\leq j \leq g$.
With these we get the adjoint operators
$\paradj=\parad{j,i}:C^{i+1}(X)\ra C^i(X)$,
and hence the partial laplacians
\[ \Lapj = \Lap{j,i} =
\parad{j,i}\pard{j,i} + \pard{j,i-1}\parad{j,i-1}:
C^i(X) \ra C^i(X) \,.\]
Each $\Lap{j,i}$ is a self-adjoint nonnegative
operator. For $i$ fixed they all commute, and one
gets a combinatorial harmonic theory (compare
\cite{St}).


When $X$ is infinite, these notions extend
to $\Ltwo$\/-cochains. When $X=\Tr$ is an
$\vr = (r_1,\dots,r_g)$\/-regular product of
trees,  Kesten's $1$\/-dimensional results 
\cite{Kes} extend and we get that each 
$\lambda$ in the spectrum of
$r_j\Id - \Lapj $ acting on $\Ltwo$\/-cochains 
of $\Tr$, satisfies
$|\lambda| \leq 2\sqrt{r_j-1}$. As in the
$1$\/-dimensional case we say that a
\rtype-regular cubical complex $\cX$ is
Ramanujan if the eigenvalues of $r_j\Id - \Lapj$
on $\cX$ are $\pm r_j$ or satisfy the same
properties for each $j$.

One justification for this definition in the
$1$\/-dimensional case is the Alon-Boppana result,
which shows that these bounds are essentially the
best possible for the trivial local system. We 
generalize this result under a 
natural hypothesis. Another parallel
with the $1$\/-dimensional case is that when $X$ 
is finite, connected,
and uniformized by a lattice $\Gamma$ in a product
$\prod_{1\leq j \leq g} G_j$ of $p$\/-adic
$\SL_2$\/'s, the Ramanujan property is equivalent
to the following condition. For a $(g-1)$\/-cube
$\sigma$ in the universal covering complex, which
is constant in the $j$\/th direction, let
$\Gamma_\sigma$ be the stabilizer
of $\sigma$ in $\Gamma$. Then no nontrivial 
unramified complementary series representations
should appear in
$\Ltwo(\Gamma_\sigma\backslash G_j)$ for any
$j$ and $\sigma$.

A fundamental problem in the $1$\/-dimensional 
case is to construct explicitly an infinite family 
of Ramanujan graphs of a fixed regularity. The 
only such families known  have
regularities of the form $k=q+1$, with $q$ 
a power of a prime. The standard
LPS examples (\cite{LPS}) depend in fact on the 
Ramanujan-Petersson conjectures for weight $2$
holomorphic cusp forms on $\GL_2(\QQ)$. These were
reduced to the Weil bounds for curves by Eichler, 
Shimura, and Igusa.  They only give examples
with $q$ a prime. The case of any prime power
was subsequently handled in \cite{Mor} via 
function fields. Ramanujan local systems were 
defined and constructed in \cite{JL6} with 
$q$ a prime. This construction requires Deligne's
results on the Ramanujan-Petersson conjecture for
cusp forms of higher weight. They depend on the Weil 
bounds for higher dimensional varieties.

In the higher dimensional case a product of 
Ramanujan graphs always gives a Ramanujan 
cubical complex, and the problem is to construct 
irreducible examples. In this work we construct 
infinite families of Ramanujan
regular cubical complexes, with fixed 
regularities $(r_1,\ldots,r_g)$, where
$r_j = q_j+1$ and the $q_j$\/'s are any
prime powers. As in the
LPS examples, our construction uses quaternion
algebras, this time over totally real number fields.
It reduces the Ramanujan property to the 
Ramanujan-Petersson  conjecture for certain 
holomorphic Hilbert modular forms. The Ramanujan-Petersson
conjecture seems not to be known for all holomorphic
Hilbert modular forms. In particular, the results of
\cite{BL} do not help us, because they establish the 
bounds outside a finite set of primes which is 
inexplicit, and a priori may depend on the form.
However, in the literature this conjecture is proved 
under additional hypotheses, which still permit us
to construct infinite families of examples with 
arbitrary prime powers for regularities.
If all the $q_j$\/'s are powers of the same prime, 
one could give, as before, an alternate construction 
using function fields. As in \cite{JL6}, we also
construct Ramanujan local systems on our examples.
Among the interesting features of our examples is
the fact that their cohomology vanishes except in
the top dimension $g$ (and $0$). In addition, 
one can bound the girth from below as in the
$1$\/-dimensional case (\cite[Theorem 7.3.7]{Lub}.

As in \cite{JL6}, this theory is valid for hermitian
local systems over cubical complexes. We similarly
give examples of Ramanujan local systems on cubical
complexes.

We gratefully acknowledge the support of the
U.S.-Israel Binational Science Foundation during
the course of this work. Additionally the first 
author was partially supported by grants from
the N.S.F. and PSC-CUNY. We thank Y. Glasner for
helpful comments on the manuscript.

\section{The harmonic theory of regular cubical
complexes}
\soussection{Cubical complexes}
In this article we will study cell complexes
locally isomorphic to a finite product of
regular trees, in which the order and the
regularities of the factors are globally constant.

Let $g$ be the dimension of such a complex
$\cX$, and let $r_j$ be the regularity of the
$j$\/th factor tree ($1\leq j\leq g$). We will 
call such a complex $\cX$ an \rtype-regular
cubical complex. Each cell in $\cX$ is a cube of
dimension $\leq g$. The $0$\/-dimensional cells
are called vertices and the $1$\/-dimensional
cells are called edges.

For each subset $I\subset\{1,\ldots,g\}$, the 
$I$\/-cubes, or cubes of direction $I$ of $\cX$,
are the products of edges from the factors $i\in I$
with vertices from the factors not in $I$. Each such cube
has $2^{|I|}$ orientations, and we denote the set
of oriented $I$\/-cubes by $\ocell{I}$. There are
{\em bottom} and {\em top} maps
$\botj,\topj:\ocell{I}\ra\ocell{I-\{j\}}$
and
{\em inversion} maps $\invj:\ocell{I}\ra\ocell{I}$
for
any $j\in I$. These are subject to the following
axioms:
\begin{enumerate}
\item $\{\invj\}_{j\in I}$ generate a group,
isomorphic to the group $\ZtwoI$ of maps from
$I$ to $\ZZ/2\ZZ$, which acts simply transitively
on the orientations of each $I$\/-cube.
\item $\topj\inv{j'}=\inv{j'}\topj$ and
$\botj\inv{j'}=\inv{j'}\botj$ for $j\neq j'$.
\item $\topj\inv{j} = \botj$ (and hence also
$\botj\inv{j} = \topj$) for all $j$.
\item Any oriented $I$\/-cube is the $j$\/th top
of precisely $r_j$ oriented $I\cup \{j\}$\/-cubes
for $j\not\in I$.
\end{enumerate}

Geometrically, these combinatorial conditions
mean that $\cX$ is locally isomorphic to the
ordered product $\Tr=\prod_j \Tr_{j}$ of $g$
regular trees of respective regularities $r_j$.
In particular, when $\cX$ is connected these
conditions hold if and only if the universal
cover $\cXtil$ is isomorphic to $\Tr$, with the
covering transformations preserving the order
of the tree factors. The if part is clear, and
the only if part holds because $\Tr$ is simply
connected (in fact, explicitly contractible),
and locally isomorphic to $X$. Hence $\Tr$ is
isomorphic to $\cXtil$ by the uniqueness of
the universal cover. The directions are
preserved under the covering map, and hence
the deck transformations form indeed a subgroup
of $\prod_j \Aut(\Tr_{j})$.

We set $\Ver=\ocell{\emptyset}$,
$\edgo{j}=\ocell{\{j\}}$, and
$\edgo{}=\sqcup_j \edgo{j}$. For an oriented edge
$e$ of direction $\{j\}$ we write
$\botj(e)=o(e)$, $\topj(e)=t(e)$, and
$\invj(e)=\ov{e}$. We will refer to these as the
the origin, the terminal vertex, and the 
oppposite edge respectively. If necessary we will
indicate the dependence of all these objects on
$\cX$.

\begin{examples}{\rm
\begin{enumerate}
\item The unit cube in $\RR^g$ is the product 
of $g$ intervals, which are \mbox{$1$\/-regular} trees.
\item The  tiling of $\RR^g$ by unit cubes with
integer vertices is the product of $g$ lines,
viewed as $2$\/-regular trees.
\item The graphs in the sense of \cite{Ser1} which
are $k$\/-regular are precisely the $k$\/-regular
cubical complexes (of dimension $1$).
\item A finite product of regular cubical
complexes is a regular cubical complex.
\end{enumerate} }
\end{examples}

A connected regular cubical complex is called
irreducible if it has no finite unramified
cover by a product of regular cubical complexes
of positive dimension.

In many important cases, the vertices of an
\rtype\/-regular cubical complex $X$ come with
{\em parities}: these are maps $p_j$ from the
vertices to $\{0,1\}$ which satisfy
$p_j(\topi e) = p_j(\boti e)$ if and only if
$i\neq j$ for any edge $e$ of direction $\{i\}$ of
$X$, $1\leq i,j \leq g$. On a complex with
parities we give each cube a canonical
orientation by agreeing that its bottom-most
vertex has all its parities equal to $0$. In
the $1$\/-dimensional case we recover the 
notion of a bipartite graph.

\soussection{Local systems}
A real/complex {\em local system}, or a
{\em flat vector bundle}, on a regular cubical
complex $\cX$ depends only on the cells of
dimension $\leq 2$. For the $1$\/-dimensional
case, see, e.g., \cite{JL6}. In the general
case, a local system $\cL$ on $\cX$ consists
of a real/complex vector space $\cL(v)$
for any vertex $v$ of $\cX$. In addition,
for any oriented edge $e$ one is given a linear
isomorphism $\cL_e:\cL(o(e))\ra\cL(t(e))$
so that $\cL_{\ov{e}} = \cL_e^{-1}$. The
$\cL_e$\/'s must also satisfy the {\em
flatness condition}\/: for any $2$\/-cell
of direction $\{j,j'\}$ in $\cX$ we require
\begin{equation}
\label{flat}
\cL_{\topj}\cL_{\botf{j'}} =
               \cL_{\topf{j'}}\cL_{\botj} \,.
\end{equation}

The local system is {\em metrized} if each
{\em fiber} $\cL(v)$ is a finite dimensional
(positive definite) inner product space and all
{\em transition maps} $\cL_e$ are isometries.
The (metrized) trivial local system $\triv{V}$,
for a finite dimensional (inner product) vector
space $V$, has  all fibers $V$ and all transition
maps  the identity. Over a contractible
space all local systems admit a trivialization,
i.e., an isomorphism to some $\triv{V}$. The 
notions of maps and direct sums make sense for
(metrized) local systems over a fixed space 
(see, e.g., \cite{JL6}). A local system
is {\em irreducible} if it is not the direct sum
of nonzero local systems. Every metrized local
system is the direct sum of irreducible ones and
a sub local system has an orthogonal complement.

Let $X_1$, $X_2$ be regular cubical complexes
with (metrized) local systems $\cL_i$ on $X_i$.
The product $X = X_1 \times X_2$ is naturally a 
regular cubical complex on which we have a
(metrized) local system
$\cL = \cL_1\boxtimes\cL_2$, called the external
product of the $\cL_i$\/'s. It is irreducible
if and only if both $\cL_i$\/'s are irreducible,
and if $X_1$ and $X_2$ are connected, then any
irreducible local system on $\cX$ is an external
 product.

Let $\cX$ be a connected regular cubical complex
and let $v_0$ be a base vertex. The universal
cover $\cXtil$ is a product of regular trees
$\{\Tr_{j}\}_{1\leq j\leq g}$, and the fundamental
group $\Pi=\pi_1(\cX,v_0)$  is discrete in
$\prod_j\Aut(\Tr_j)$ which acts properly on $\cXtil$
with quotient $\cX$. (Metrized) local systems $\cL$
on $\cX$ are equivalent to (orthogonal or unitary)
representations $\rho_\cL$ of $\Pi$ on $\cL(v_0)$.
We reconstruct $\cL$ as
$\Pi\backslash(\cXtil\times V)$.

A similar construction is possible in the
disconnected case: if $\Pi$ acts freely on
a regular cubical complex $\cXtil$ preserving
directions, and $\rho$ is an (orthogonal or
unitary) representation of $\Pi$ on $V$, then
the quotient $\cX=\Pi\backslash\cXtil$ is
again a regular cubical complex and we get a
(metrized) local system
$\cL=\Pi\backslash(\cXtil\times V)$. One gets
this way all local systems on $\cX$ whose
pullback to $\cXtil$ is trivial.

\soussection{Cochains}
Let $\cX$ be an \rtype-regular cubical complex
and let $\cL$ be a local system on $\cX$. For
a subset $I\subset \{1,\ldots,g\}$ and
$\sig\in\ocell{I}$ denote by $o(\sig)$ the
iterated bottom vertex 
$(\prod_{j\in I}\botj)(\sig)$, and for $j\in I$ 
let $e_j(\sig)$ be the oriented edge of direction $\{j\}$
defined by
$e_j(\sig) = (\prod_{j'\neq j}\botf{j'})(\sig)$.
By the combinatorial conditions these
are well defined, i.e., do not depend on the order
of the operations. Now set $\cL(\sig)=\cL(o(\sig))$
and $\cL_{\sig,j} = \cL_{e_j(\sig)}$. By
definition, the space of $I$\/-cochains $C^{I}(\cX,\cL)$
is the space of completely alternating collections
\[s=\set{s(\sig)}_{\sig\in\ocell{I}}
\in \prod_{\sig\in\ocell{I}}\cL(\sig)\,.\]
The condition of being completely alternating is
vacuous for $0$\/-cochains. Otherwise, a collection
$s$ as above is completely alternating if for any
$j\in I$ we have
$s(\invj(\sig)) = - \cL_{\sig,j}s(\sig)$.
For $\cL = \cT_V$ we can view the $I$\/-cochains
as maps from $\ocell{I}$ to $V$.

For an integer $0\leq i \leq g$  put
$C^i(\cX,\cL)= \oplus_{|I|=i} C^I(\cX,\cL)$.
The partial boundary operators
$\parj = \pard{j,i}: C^i(\cX,\cL)
                   \ra C^{i+1}(\cX,\cL)$
are defined by linearity and the formula
\[\parj s(\sig) = \begin{cases}
\cL_{j,\sig}^{-1}s(\topj(\sig)) - s(\botj(\sig)) &
\text{for $\sig\in\ocell{I}$ with
                 $|I|=i+1$ and $j\in I$} \\
0 & \text{for $\sig\in\ocell{I}$ with
           $|I|=i+1$ and $j\not\in I$}
\end{cases} \]
for any $s\in C^i(\cX,\cL)$. These are well defined
and satisfy $\parj\pard{j'} = \pard{j'}\parj$
for all $j,j'$. Clearly $\parj^2=0$. We will denote
by $\pard{j,I}$ the restriction of $\pard{j}$ to
$C^I(\cX,\cL)$.
If $j\in I$, then $\pard{j,I}=0$.
 As in
\cite[II.1]{Ser}, the total boundary of
$s\in C^I(\cX,\cL)$ is defined by 
$\der(s) = \der_I(s) = 
    \sum_{j\not\in I} (-1)^{\alpha_I(j)}\pard{j,I}$,
where $\alpha_I(j)$ is the  place of $j$ in
$I\cup \{j\}$. As
usual $\der^2=0$. We get the spaces of
cohomology with respect to $\parj$ and $\der$
by the usual formulas
$H^i_j(\cX,\cL) = \Ker\pard{j,i}/\Ima\pard{j,i-1}$
and $H^i(\cX,\cL) = \Ker\der_i/\Ima \der_{i-1}$. 
We can moreover break $H^i_j(\cX,\cL)$ as a direct
sum of subspaces $H^I_j(\cX,\cL)$ with $|I| = i$
defined by restricting $d_j$ to $C^I(\cX,\cL)$.
When $\cL$ is
the trivial local system $\cT_\RR$, it is usually
omitted from the notation.

The operators
$\paradj = \parad{j,i}:
        C^{i+1}(\cX,\cL) \ra C^i(\cX,\cL)$
are defined by
\[ \parad{j,i}s(\sig) =
         \sum_{\topj(\tau) = \sig}
                \cL_{\tau,j}^{-1} s(\tau)\,.\]
These operators satisfy the analogous relations
$\paradj\parad{j'} = \parad{j'}\paradj$ and 
$(\paradj )^{2}=0$. We let $\parad{j,I}$
denote the restriction of $\parad{j,i}$ to
$C^{I\cup\{j\}}(\cX,\cL)$ for $j\not\in I$. Finally, for
$t \in C^I(\cX,\cL)$ we define 
$\derad(t) =
     \sum_{j\in I} (-1)^
     {\alpha_{I-\{j\}}(j)}\parad{j,I-\{j\}}$.
Again $(\derad)^2=0$.

\soussection{Hodge Theory}
From now on we consider only metrized local
systems.  Let $C_2^I(\cX,\cL)$ denote the Hilbert
space completion of the space of $I$\/-cochains 
with coefficients in $\cL$, with respect to the
pre-Hilbert norm
\[\|s\|^2 = 2^{-|I|}
     \sum_{\sig\in\ocell{I}} \|s(\sig)\|^2\,.\]
We let $C_2^i(\cX,\cL)$ be the orthogonal sum of
the corresponding $C^I$\/'s. It is clear that the
operators previously defined (e.g.,
$\pard{j,i}$, $\derad_i$,  etc.) induce bounded
operators, denoted by the same letters, on
the corresponding spaces $C_2^*(\cX,\cL)$. 
In particular we have $\|\pard{j,i}\|^2 \leq r_j$
and $\|\parad{j,i}\|^2 \leq r_j$.
To see this, notice that for $|I|=i$ and 
$s\in C^I(\cX,\cL)$ we
have $\pard{j,i}s = 0$ if $j$ is in $I$. Otherwise
set $I'=I\cup \{j\}$. Then
 \begin{align*}
\|\pard{j,i} s\|^2 & =
 2^{-i-1}\sum_{\sigma\in\ocell{I'}}
 \| \cL_{j,\sigma}^{-1}s(\topj(\sigma)) -
                            s(\botj(\sigma)) \|^2\\
 & \leq  2^{-i+1}\sum_{\sigma\in\ocell{I'}}
 \| s(\botj(\sigma)) \|^2 
  =   r_j 2^{-i}\sum_{\tau\in\ocell{I}}
 \| s(\tau) \|^2\\
 & = r_j \|s\|^2
\end{align*}
(compare \cite[Lemma 1.2]{JL6}). This implies
the bound
$\|\pard{j,I}\|^2 \leq r_j$.
Since the $C^I(\cX,\cL)$\/'s are mutually 
orthogonal, the same bound holds for the
direct sum operator $\pard{j,i}$.
The case of $\parad{j,i}$ is 
similar and we omit it.

It is routine and easy to verify that 
$\parad{j,I}$ is the adjoint of $\pard{j,I}$.

The laplacians
$\Lap{j,i} : C^i(\cX,\cL) \ra C^i(\cX,\cL)$
are defined as usual by
$\Lap{j,i} = \pard{j,i-1}\parad{j,i-1} +
            \parad{j,i}\pard{j,i}$.
On each $C^I(\cX,\cL)$ this simplifies
to $\pard{j,I-\{j\}}\parad{j,I-\{j\}}$ if 
$j\in I$
and to $\parad{j, I}\pard{j,I}$ if $j\not\in I$.
These are commuting bounded self-adjoint 
nonnegative operators.
The Laplacian $\Lapt{i}$
is defined by
\[
\Lapt{i}=\Lap{1,i}+\Lap{2,i} + \cdots +\Lap{g,i} =
\der_{i-1}\derad_{i-1} +\derad_{i}\der_{i}\  .
\]
The restriction of $\Lap{j,i}$ 
($\Lapt{i}$) to
$C^I(\cX,\cL)$ will be denoted by $\Lap{j,I}$
($\Lapt{I}$).
All these are also bounded, self-adjoint, 
and nonnegative. By definition, the space of
$\bul$\/-harmonic forms for $\bul = i$ or $I$,
$\cH^\bul(\cX,\cL)\subset C^\bul(\cX,\cL)$, is
$\Ker \Lapt{\bul}$. We likewise set
$\cH^{\bul}_{j}(\cX,\cL) = \Ker \Lap{j,\bul}$. We
now have the following routine
\begin{proposition}
Let $\cX$ be a finite regular cubical complex, 
and let $\cL$ be a metrized local system on 
$\cX$. Then
\begin{enumerate}
\item $\cH^i(\cX,\cL) = 
\Ker\,\der_i\cap\Ker\,\derad_{i-1}$ and
$\cH_{j}^{i}(\cX,\cL) =
\Ker\,\pard{j,i}\cap\Ker\,\parad{j,i-1}$.
Moreover
$ \cH_{j}^{I}(\cX,\cL) =
\Ker\,\pard{j,I}$
if $j\notin I$ and  $ \cH_{j}^{I}(\cX,\cL) =
\Ker\,\parad{j,I-\{j\}}$ if $j\in I$.
\item We have orthogonal sum decompositions 
{\rm (}the {\em Hodge Decomposition}{\rm )}
\begin{align*}
C_2^i(\cX, \cL) & = \cH^i(\cX, \cL) \oplus \Ima\,\der_{i-1}
         \oplus \Ima \,\derad_{i}
         \quad , \\
C_{2}^{i}(\cX, \cL) & = \cH_{j}^{i}(\cX,\cL) \oplus
        \Ima\,\pard{j,i-1} \oplus
                  \Ima \,\parad{j,i}\ ,\quad
                  \text{and}\\
C_{2}^{I}(\cX, \cL) & = \cH_{j}^{I}(\cX,\cL) \oplus
                  \Ima \,\parad{j,I} \quad\text
                  {for $j\not\in I$}\\
                  & =  \cH_{j}^{I}(\cX,\cL) \oplus
                   \Ima \,\pard{j,I-\{j\}} \quad\text
                  {for $j\in I$}\ .\\
\end{align*}
\end{enumerate}
\end{proposition}

\section{The spectrum of the laplacian}
Our main interest will be in the eigenvalues 
of the $\Lap{j,i}$\/'s. First, we have 
the following
\begin{lemma}
For every subset $I\subset \set{1,\dots,g}$ the
$I$\/-cochains $C^I(X,\cL)$ are preserved by
$\Lap{j,i}$ for all $j$, with $i = |I|$. Moreover,
for $j\not\in I$ the maps $\parj$ and $\paradj$
induce isomorphisms between the
$\lambda$\/-eigenspace of $\Lapj$ on $C^I(X,\cL)$
and on $C^{I\cup \set{j}}(X,\cL)$ for any
$\lambda \neq 0$.
\end{lemma}
\begin{pf}
Suppose $\Lapj c = \lambda c$ and
for $\Lapj c' = \lambda c'$ for
$c\in C^I(X,\cL)$ and $c'\in C^{I\cup j}(X,\cL)$,
with $j\not \in I$. Then
$\Lapj \parj c = \parj \paradj \parj c =
 \lambda \parj c$ and likewise $\Lapj \paradj c'
= \lambda \paradj c'$. Since
$\paradj \parj $ and
$\parj \paradj $ are both multiplication by
$\lambda\neq 0$ on all such $c$, $c'$ respectively, 
the lemma follows.
\end{pf}
\noindent It follows that the study of the nonzero
eigenvalues and eigenspaces of $\Lapj$ on any
$C^I(X,\cL)$ can be reduced to 
the case $j\not \in I$. In fact, it is also possible
to reduce the study of the full diagram
\[C^I(\cX,\cL)  \overset{\parj}
     {\underset{\paradj}{\rightleftarrows}}
  C^{I\cup\{j\}}(\cX,\cL)\,,\]
to the case when $X$ is an $r_j$\/-regular graph
and $I$ is empty. To do this, we make for convenience 
the following hypothesis.
\[ \mbox{PAR:\quad\quad $X$ has parities $p_j$ for all
   $1\leq j \leq g$.} \]
Under hypothesis PAR we define, for any $j\not\in I$,
a graph $\Gr_j(X) = \Gr_{j,I}(X)$ and a (metrized)
local system $\cL = \cL_{j,I}(X)$
as follows. The vertices of $\Gr_j(X)$ are the cubes of type
$I$ of $X$; its (oriented) edges are the (oriented)
cubes of type $I\cup \{j\}$,  (with the $I$ part of the
orientation all $0$\/'s).
For an $I$\/-cube $\sigma$ of $X$
define $\cL_{j,I}(X)(\sigma) = \cL(v_\sigma)$,
where $v_\sigma$ is the bottom-most vertex of $\sigma$.
The transition maps are then those induced from $\cL$.
The compatibilities of the corresponding operators
$\parj$, $\paradj$, $\Lapj$ with $\pard{}$,
$\parad{}$, and $\Lap{}$ are obvious from the definitions.
In any mention of these graphs hypothesis PAR will be assumed.

Let $\Gr$ be a a locally finite graph and $\cL$\
be a local system on $\Gr$.  The star operator
$S(\Gr, \cL): C^{0}(\Gr, \cL)\ra C^{0}(\Gr, \cL)$
is defined to be the map sending a cochain
\mbox{$s\in C^{0}(\Gr ,\cL)$} to the cochain
$S(\Gr, \cL)(s)$ given by
\[
S(\Gr, \cL)(s)(v)=\sum_{\{e\,|\, t(e)=v\}}
\cL_{e}^{-1}s(o(e)) \ .\]
In case the graph $\Gr = \Gr_{j,I}(X)$ we shall
call $S_{j,I}=S(\Gr_{j,I}(X),\cL)$ the $j$\/th star
operator 
and put $S_{j,i}=\oplus_{|I|=i}S_{j,I}$. The 
operator $S_{j,I}$ then acts on
the space $C^I(\cX,\cL)=C^{0}(\Gr_{j,I}(X),\cL)$
and accordingly $S_{j,i}$ acts on the space
$C^{i}(\cX,\cL)$.
Clearly the norm of $S_{j,I}$ is bounded
by $r_j$. 
Since $\Lap{j,I}$ can be viewed as $r_j - S_{j,I}$
(compare \cite[Proposition 2.20]{JL6}),
we get the following
\begin{proposition}
The norm of $\Lap{j,I}$ restricted to $C^I(\cX,\cL)$
is bounded by $2r_j$ and
its spectrum is contained in $[0,2r_j]$.
\end{proposition}
\noindent As a corollary, the same bound is valid for the norm
of $\Lap{j,i}$ (use the orthogonal decomposition of
$C^i(\cX,\cL)$ to the $C^I(\cX,\cL)$\/'s, treating
separately the cases $j\in I$ and $j\not\in I$).

In dimension $1$, the classical result of Alon and Boppana
(see, e.g., \cite[Proposition 4.5.4]{Lub}) shows
that the nontrivial eigenvalues of the star operator for
the trivial local system are
essentially bounded in terms
of the norm of the star operator on the universal
covering regular tree (see also \cite{Kes} or 
\cite{Bu}). This result is the following:
\begin{proposition}{\rm (Alon-Bopanna)}
Let $\Gr_n$ be a family of finite 
$r$\/-regular connected graphs whose number of 
vertices goes to $\infty$ with $n$. Let 
$L_{2,0}(\Gr_n)$ denote the space of zero-sum
$\Ltwo$ zero cochains. Then 
\[\liminf_{n\ra \infty} 
\| S(\Gr_n)|_{L_{2,0}(\Gr_n)} \|
  \geq 2\sqrt{r-1}\,.\]
\end{proposition}
In the higher dimensional case set $\bul = i$
or $I$, and define for $j\not\in I$
\[\mu_{j,\bul}(X,\cL) = \max\set{|\mu|\,|
\,\, \mu \text{ is an eigenvalue of } S_{j,\bul}
\text{ different from }\pm r_j }\,. \]
Here $\cL$ is any metrized local system on $X$.
Put also  
\[\mu_{j,\bul}(X) = \mu_{j,\bul}(X,\triv)\,.\]
In particular, $\mu_{j,\bul}(X,\cL)$ lies in
the interval $[-r_j,r_j]$.
The Alon-Boppana result now implies the
following $g$\/-dimensional version:
\begin{proposition}
Let $X_n$ be a sequence of $(r_1,\dots,r_g)$\/-regular
connected cubical complexes.  If the  
number of vertices of each connected component 
of $\Gr_{j,I}(X_n)$ tend to $\infty$ with $n$, then 
$\liminf_{n} \mu_{j,I}(X_n) \geq 2\sqrt{r_j-1}$.
\end{proposition}

In analogy with the $1$\/-dimensional case
(\cite[Definition 2.25]{JL6}) we make the following
\begin{definition}
A local system $\cL$ on a regular cubical 
complex $X$ is Ramanujan if
$\mu_{j,I}(X,\cL) \leq 2\sqrt{r_j-1}$ for all 
$j$ and $I$ as above.

We say that $X$ is Ramanujan if the trivial local 
system $\triv{\CC}$ on $X$ is Ramanujan.
\end{definition}
\noindent It is of course possible to define the notion of
an $(j,I)$\/-Ramanujan system, so that $\cL$ is
Ramanujan if and only if it is $(j,I)$\/-Ramanujan
for all $(j,I)$ ($j\not\in I$).

Many examples are known in the $1$\/-dimensional 
case (see, e.g., \cite{Lub}). Since the external 
tensor product of Ramanujan local systems
is Ramanujan, we get many examples in the
higher-dimensional case as products, or more
generally by pulling back a product system 
from a finite unramified cover. We will say
that a connected regular cubical complex is
reducible if a finite unramified cover of it is a
product of cubical complexes of lower dimensions.
Otherwise we will say that $X$ is irreducible.

When $X$ be connected, we can write 
$X$ as a quotient $X = \Gamma\backslash \Tr$
of a product $\Tr = \prod_j \Tr_j$ of trees.
Then $X$ is  
{\em reducible} if and only if the following
two conditions are satisfied: a) $\Tr$ is a
product $\Tr = \Tr' \times \Tr''$, with $\Tr'$, 
$\Tr''$
products over complementary subsets of the 
$j$\/'s; and b) A subgroup of finite index of
$\Gamma$ is compatibly a product
$\Gamma' \times \Gamma''$, so that the
corresponding finite unramified cover of $X$ is
$\Gamma'\backslash \Tr' \times
                     \Gamma''\backslash \Tr''$.

A metrized local system is always an orthogonal 
sum of irreducible local systems. The irreducible 
(metrized) local systems 
over a product (connected) complex are precisely
the external tensor products of irreducible ones. 
The challenge is then to find Ramanujan local
systems over irreducible cubical complexes.
Even when a graph or a cubical complex is Ramanujan,
there are local systems on it which are not: a
generic deformation of the transition maps of the 
trivial local system is an example.

Let $F_1$,\dots, $F_g$ be nonarchimedean local
fields, and set
$G_j = \PGL_2(F_j)$ and $G = \prod_j G_j$.
Let $G_+$ be the elements of $G$ whose $j$\/th
components have determinants with even valuations
for each $j$.
 The standard trees $\Tr_j$ associated to the
$G_j$\/'s are 
$r_j = (q_j + 1)$\/-regular, with $q_j$
the cardinality of the residue field of $F_j$.
Set $\Tr = \prod_j \Tr_j$.
For a discrete subgroup $\Gamma$ of $G_+$ the 
quotient complex $\Gamma\backslash \Tr$ satisfies
hypothesis PAR. For any
\mbox{$I\subset \set{1,\dots,g}$}, $G$ acts on
$\Tr_I =\prod_{j\in I} \Tr_j$ and on
$\Tr^I =\prod_{j\not\in I} \Tr_j$, and
$\Tr = \Tr_I \times \Tr^I$.
For an $I$\/-cube $\sigma\in\Sigma^I(\Tr_I)$ let
$\Gamma_\sigma$ denote the projection to 
$G^I = \prod_{j\not\in I} G_j$ of the stabilizer
of $\sigma$. We now have the following
\begin{proposition}
\label{compl}
{\rm 1.}\ Let $\Gamma$ be a discrete, cocompact,
torsion-free subgroup of $G$. Let $\cL$ be the
{\rm (}metrized\/{\rm )} local system $\cL$ on
$X = \Gamma \backslash \Tr$ corresponding
to a unitary representation $\rho$ of $\Gamma$ on
a {\rm (}finite-dimensional\/{\rm )} space $V$.
Then $\cL$ is $(j,I)$\/-Ramanujan for $j\not\in I$ 
if and only if the following condition holds. For 
any $I$\/-cube $\sigma$ of $\Tr^{\set{j}}$, no
{\rm (}nontrivial\/{\rm )} representations of $G$
of the unramified complementary series appear in
$\Ltwo (\Gamma_\sigma \backslash G_j\times V)$, 
where $G_j$ acts through its right action on the 
$G_j$\/-factor. \\[0.05in]
{\rm 2.}\ The local system $\cL$ is Ramanujan on $X$ if
and only if no nontrivial representations of $G$
of the unramified complementary series appear in 
$\Ltwo (\Gamma_\sigma \backslash G_j\times V)$, 
for any cube $\sigma$ of $\Tr^{\set{j}}$ of direction
$I^j = \set{1,\dots , g}- \set{j}$.

\end{proposition}
\begin{pf}
1.\ After the choice of a maximal unramified 
compact
subgroup of $G$, the unramified vectors in
$\Ltwo (\Gamma_\sigma \backslash G_j\times V)$
can be identified with
$\Ltwo (\Gamma_\sigma \backslash 
                     \Ver \Tr_j \times V)$,
and the action of the $j$\/th factor 
Hecke operator $T_{v_{j}}$
corresponds to the action of the $j$\/th star
operator (compare Step 1 in the proof of
\cite[Theorem 3.4]{JL6}). Our claim now
follows from Satake's reformulation of the
Ramanujan-Petersson conjecture (\cite{Sat}).\\[0.05in]
2.\ For $j\not\in I\subset \set{1,\dots,g}$, an
$I^j$\/-cube has $2^{g-1-|I|}$ faces of direction
$I$. Hence the diagram 
\[C^I(\cX,\cL) \overset{\parj}
     {\underset{\paradj}{\rightleftarrows}}
  C^{I\cup\{j\}}(\cX,\cL)\,,\]
embeds (in $2^{g-1-|I|}$ ways) into the corresponding
diagram with $I$ replaced by $I^j$. Hence to
verify the Ramanujan property it suffices to verify
the $(j,I^j)$\/-Ramanujan property for all 
$1\leq j\leq g$. For any $I$ as above we have
the decomposition into connected components
\[\Gr_{j,I}(\cX,\cL) = 
   \sqcup \Gamma_\sigma 
        \backslash \Tr^{I\cup\set{j}} \,,\]
the union taken over representatives $\sigma$ 
of $\Gamma \backslash \Sigma^I(\Tr^{\set{j}})$. It
suffices therefore to verify the Ramanujan
property for $\cL$ over each of these components
for $I = I^j$, and our claim follows from the first 
part.
\end{pf}

\section{Ramanujan local systems arising from
quaternion algebras}
\label{qua}
Let $F$ be a totally real field of degree $d$ 
over $\QQ$ with ring of integers $\cO_F$. Let
$\infty_i$, $i = 1,\dots,d$, be the infinite 
places of $F$. For an algebraic group $\bH$
over $F$ we let $H_F$ be the $F$\/-rational
points of $\bH$, let $H_\AA$ be the adelic
points of $\bH$, and let $H^\rf$ denote the
finite-adelic points of $\bH$. For any
finite set of finite places $v_1,\dots,v_n$
of $F$ let $H^{\rf,v_1,\dots,v_n}$
denote the finite-adelic points of $\bH$
without the ${v_1,\dots,v_n}$ component.

Let $B$ be a totally definite 
quaternion algebra over $F$ with reduced norm
$\Nm = \Nm_{B/F}$. (A general reference for
quaternion algebras over totally real fields is
\cite{Vig}.) Let $\bG$ be the
algebraic group over $F$ associated with the 
multiplicative group $B^\times$. Let $q_v$
denote the cardinality of the residue field of
the completion $F_v$ of $F$ at a nonarchimedean
place $v$.

With $\HH$ denoting the Hamilton quaternions,
let $\Symm^k$ denote the $k$\/th symmetric
power of the $2$\/-dimensional representation 
of $\HH^\times$ obtained by identifying
$\HH\otimes_\RR\CC$ with $\Matwo(\CC)$. We have
$G_\infty := \prod_i \bG(F_{\infty_i})
                      \simeq (\HH^\times)^d$.
The irreducible representations of $G_{\infty}$
are then
$\Symm^{\vk,\vs} = \otimes_i
     (\Nm_{B_{\infty_i}/\RR})^{s_i} \Symm^{k_i}$,
with $\vk = (k_1,\dots,k_d)$ and
$\vs = (s_1,\dots,s_d)$, the $s_i$ being any complex
numbers. We denote by $V^{\vk,\vs}$ the space on which
$\Symm^{\vk,\vs}$ acts. They can be unitarized if and
only if the center of $G_\infty$ acts through a 
unitary character, and this happens if and only 
if $2 s_i + k_i$ is imaginary for each $i$. 
In this event
there is a $G_{\infty}$-invariant positive 
definite hermitian inner
product on $V^{\vk ,\vs}$ which is unique up to 
a scalar.  In what
follows we shall also assume that $\vs = -\vk/2$.
Then an element $(z_1,\dots,z_d)$ in the center
of $G_\infty$ acts via $\prod_i \sign(z_i)^{k_i}$.
In particular, an element $x$ of the center 
$F^\times$ of $G_F \subset G_\infty$ acts via
$\sign\bigl(\Nm_{F/\QQ}(x)^k\bigr)$ if all 
the $k_i$\/'s have the same parity $k$.

Let $\bZ$ be the center of $\bG$,
viewed as an algebraic group over $F$. Let $v_j$,
$j=1,\dots, g$, be distinct nonarchimedean
places of $F$ where $B$ is unramified, and fix
identifications 
\begin{equation}
\label{iden}
B_{v_j} \simeq \Matwo(F_{v_j})\,.
\end{equation}
Set 
$G_j = \GL_2(F_{v_j})$, $G=\prod_{j}G_j$,
and $K_j = Z_{v_j}\GL_2(\cO_{F,v_j})$.
We view the vertices of the 
($r_j=q_{v_j}+1$\/-regular) tree $\Tr_j$ of
$G_j$ as $\Ver \Tr_j=G_j /K_j$, and set
$\Tr = \prod_j \Tr_j$.
For a compact open subgroup $K^0$ of
$G^{\rf,v_1,\dots,v_g}$
set $\Gamma = \Gamma(K^0) = G_F \cap K^0$. Then
$\Gamma$ divided by its center acts discretely on
$\Tr$, and if $K^0$ is sufficiently small the action
is free. We can therefore use the notation and
results of the previous section. Recall that the 
Ramanujan-Petersson conjecture asserts that the 
eigenvalues of the Hecke operator $T_{v_j}$
are bounded by $2\sqrt{q_{v_j}}$ on the
automorphic forms on $\GL_2(F)$ which 
are unramified at $v_j$. We now have the following
\begin{theorem}
\label{Ram}
\begin{enumerate}
\item
The cubical complex $\cX=\Gamma\backslash \Tr$
is irreducible, and the graph $\Gr_{j,I}(\cX)$
is connected for any $j\not\in I\subseteq
 \{1,\dots, g\}$.
\item Suppose that every element in $\Gamma$
which fixes a point on $\Tr$ acts trivially on
$V^{\vk,-\vk/2}$. Then the formula
$\cL = \cL^{\vk,-\vk/2} =
   \Gamma \backslash (\Tr \times V^{\vk,-\vk/2})$
defines a local system on
$X = \Gamma \backslash \Tr$, which is metrized,
irreducible, and nontrivial unless all the $k_j$
are $0$.
\item The cohomology groups $H^\ast(\cX,\cL)$
vanish except in dimensions $0$ and $g$.
\item When $\cL$ defines a local system on 
$\cX$, it is Ramanujan
provided the Ramanujan-Petersson conjecture holds
for {\em holomorphic} automorphic forms on $\GL_2(F)$
of weight $(k_1 + 2,\dots,k_g + 2)$.
\end{enumerate}
\end{theorem}
\begin{pf}
1.\ By the property of strong approximation the
projection of $\Gamma$ to each strictly partial
$G_I$ is dense, and the
$G$\/-action on the (unoriented) $I$\/-cubes
of $\Tr^{\set{j}}$ is transitive. The first fact
implies that $\cX$ is irreducible, and the second 
that $\Gr_{j,I}(\cX)$ is connected.  

\vspace{0.05in}
\noindent
2.\ The given condition clearly implies that
$\cL$ as defined is a local system on $X$.
(In our case it simply means that if $z$ is
a central element in $\Gamma$, then
$\prod_{1\leq i\leq g}
        (\sign_{\infty_i} x)^{k_i} = 1$\/.)
in our case
Since $\vs = -\vk/2$, $\rho^{\vk,-\vk/2}$ is
unitary, and hence $\cL$ is metrized. The
Eichler-Kneser strong approximation theorem
implies that the image of $\Gamma$ in $G_\infty$
contains a Zariski dense subgroup of the
group $(G_1)_\infty \simeq \SU(2)^d$ of the norm
$1$ elements of $G_\infty$, implying
the irreducibility of $\cL^{\vk,-\vk/2}$. The
nontriviality statement is clear.

\vspace{0.05in}
\noindent
3.\ This is a consequence of representation-theoretic
results of Garland and Casselman, see
\cite[Chapter 13, Proposition 3.6(i)]{BW}.

\vspace{0.05in}
\noindent
4. Hypothesis PAR holds since the elements of
$\Gamma$ have norm $1$.
We will deduce the Ramanujan property from bounds on
Hecke eigenvalues as in \cite[Theorem 3.4.2]{JL6}.
Define a finite $\vr$\/-regular cubical complex
$Y = Y(K^0)$ and a local system $\cL_Y$ on $Y$ by
\begin{equation}
\label{adel}
Y = G_F \backslash
   (\Tr \times (G^{\rf,v_1,\dots,v_g}/K^0))
   \quad\text{and}\quad
\cL_Y = G_F \backslash (\Tr \times
           (G^{\rf,v_1,\dots,v_g}/K^0)
                   \times V^{\vk,-\vk/2})\,,
\end{equation}
with the diagonal $G_F$ action.
The connected components of $Y$ are
in bijection with the id\`ele class group
$F^\times \backslash \AA_F^{\times,\rf}/
                           \Nm_{B/F}(KK^0)$
by the Eichler-Kneser strong approximation
theorem.
Moreover, the natural inclusion of $\Tr$ into
$\Tr \times G^{\rf,v_1,\dots,v_g}$ exhibits $X$
as a connected component of $Y$, with $\cL$ the
restriction of $\cL_Y$ to $X$. By
Proposition~\ref{compl}(2), the Ramanujan property
is equivalent to bounds for the ``interesting''
eigenvalues of $S_{j,I^j}$ on each
\mbox{
$C^0 (\Gr:=\Gr_{j,I^j}(Y),\cL_{Y})$, $1\leq j \leq g$}.

First observe that the stabilizer $\Iw^j$ of a fixed
$I^j$\/-cube
$\sigma$ of $\Tr^{\set{j}}$ in $G$ is the product
$\prod_{i\in I^j} \Iw_i$ of the stabilizers $\Iw_i$
of each $i$\/th edge factor of $\sigma$. Set
$K = K^0 K_j \Iw^j$. Then we are reduced to
studying the eigenvalues of the star operator
$S_{j,I^j}$ on
\begin{eqnarray*}
C^0(\Gr,\cL_Y) & = & G_F \backslash (\Ver(\Tr_j) \times
           (G^{\rf,v_j}/K^0\Iw^j)
                   \times V^{\vk,-\vk/2}) \\
    & = & G_F \backslash ((G^\rf/K^0 K_j \Iw^j)
                        \times V^{\vk,-\vk/2}) \ .
\end{eqnarray*}

Let $\cB^{\vk,-\vk/2}(\bG)$ be the space of continuous
maps $\phi: G^\rf \ra V^{\vk,-\vk/2}$ satisfying
$\phi(g^\rf x) = \rho^{\vk,-\vk/2}(g_\infty)\phi(x)$,
for any
$g = g_\infty g^\rf \in
         (G_\infty \times G^\rf) \cap G_F$,
with $G^\rf$ acting on it through right translations.
We may view $C^0(\Gr,\cL_Y)$ as the $K$\/-invariants
of $\cB^{\vk,-\vk/2}(\bG)$, with $S_{j,I^j}$ corresponding
to the Hecke operator $T_{v_j}$. We shall express
this space as a space of automorphic
forms for $B^\times$. First we decompose
under the action of $Z^\rf$:
\[
\cB^{\vk,-\vk/2}(\bG)^K =
   \oplus_{\omega \in \Omega(\bG,K^0)}
           \cB^{\vk,-\vk/2}(\bG,\omega)^K\,,
\]
with the sum over the set $\Omega(\bG,K^0)$
characters of $Z_\AA$ trivial on $Z_F$,
on the connected component of the the identity
of $Z_\infty$, and on $K \cap Z^\rf$. In
particular, it is a finite set.

As in \cite[Section 1.1]{JL5}, each
$\cB^{\vk,-\vk/2}(\bG,\omega)$ is closely related
to the space $\cA(\bG,\omega)$ of automorphic forms
on $B^\times$ in the sense of Jacquet-Langlands
\cite[Chapter 14]{JaLa}. (By the compactness of
$G_\infty^1$ these can be viewed as the
complex-valued continuous functions
$f: G_F \backslash G_\AA \ra \CC$ which are
right $G_\AA$\/-finite and which satisfy
$f(gz) = \omega(z)f(g)$ for any $g\in G_\AA$ and
$z\in Z_\AA$.) For an irreducible representation
$\rho$ of $G_\infty^1$ which is isomorphic to
$\rho^{\vk,-\vk/2}_{|G_\infty^1}$, the Peter-Weyl
theory furnishes an isomorphism
\[ \cB^{\vk,-\vk/2}(\bG,\omega) \simeq
\Hom_{G_\infty^1}(V^{\vk,-\vk/2},
           \cA(\bG,\omega)) \simeq
\Hom_{G_\infty}(V^{\vk,-\vk/2},
                  \cA(\bG,\omega))\,. \]
The second isomorphism holds because both
$\omega$ and $\rho^{\vk,-\vk/2}$ are trivial
on the connected component of
$Z_\infty \simeq \prod_i F^\times_{\infty_i}$.
The spaces involved are nonzero if and only if
the central
character of $\rho^{k,s}$ agrees with $\omega$
on $Z_\infty \simeq \prod_i F^\times_{\infty_i}$.
In the opposite direction, we can identify the
$\rho$\/-isotypical part $\cA(\bG,\omega)^\rho$
of $\cA(\bG,\omega)$ with
$V^{\vk,-\vk/2} \otimes_\CC
               \cB^{\vk,-\vk/2}(\bG,\omega)$; see,
e.g., \cite[Section 1.1]{JL5} for explicit formulas
for these isomorphisms (in the case $F = \QQ$, but
the generalization is immediate). These isomorphisms
are $G^\rf$\/-equivariant. We need
therefore to study the eigenvalues of $T_{v_j}$
on $\cA(\bG,\omega)^K$ for each
$\omega\in \Omega(\bG,K^0)$.

Let $\cA(\bG,\omega)^K_1$ be the space of
functions in $\cA(\bG,\omega)^K$ factorizing
through the norm, and $\cA(\bG,\omega)^K_2$
its orthogonal complement. Then
\[\cA(\bG,\omega)^K =
     \cA(\bG,\omega)^K_1
    \oplus \cA(\bG,\omega)^K_2\,.\]
Moreover, $\cA(\bG,\omega)^K_1$ is the sum of character
spaces $V(\chi)$, each $1$\/-dimensional, with
$\chi$ going over the id\`ele class characters
trivial on $\Nm_{B/F}(K)$ and satisfying
$\chi^2 = \omega$. On the corresponding space
$T_{v_j}$ acts as the scalar
$\lambda_\chi = (q_{v_j} + 1)\chi(\pi_{v_j})$,
where $\pi_{v_j}$ is a uniformizer for $v_j$.
Hence $\lambda_\chi^2 = r_j^2\omega(\pi_{v_j})$,
and since $\pi_{v_j}$ is in $Z_{v_j}$ we get
$\lambda_\chi^2 = \pm r_j$.

The elements of $\cA(\bG,\omega)^K_2$ lift
Hecke-equivariantly to cusp forms on $\GL_2(F)$
by the Eichler-Shimizu-Jacquet-Langlands theory.
The lift is
injective, and its image is characterized by
square integrability at each prime $v$ which is
ramified for $B$. When $v$ is
infinite, the square integrability is essentially
the same as holomorphy. Moreover, the weights
correspond as indicated: to $\Symm^k$ at an
infinite place correspond forms of weight $k+2$
at that place on $\GL_2(F)$. Hence we have
reduced the Ramanujan property for $X$ to the
Ramanujan-Petersson conjecture for holomorphic
forms on $\GL_2(F)$ of the type claimed.
\end{pf}
\begin{remarks}
 {\rm
1.\ We need only special cases
of the holomorphic Ramanujan-Petersson
conjecture. Let $B/F$; $v_1,\dots,v_g$;
$K^0$; and $\vk$ be as above. To
prove that $\cL^{\vk,-\vk/2}$ is a Ramanujan local
system on $X$ only requires that the eigenvalues
of $T_{v_1},\dots,T_{v_g}$ on those automorphic
forms on $\GL_2(F)$ which are lifts of
$\rho^{\vk,-\vk/2}$\/-isotypical forms on
$\cA(\bG,\omega)^K$ with
$\omega\in \Omega(\bG,K^0)$ satisfy the
Ramanujan-Petersson bound. These are the
automorphic forms on $\GL_2(F)$, of level
dividing $K^0$, of $\infty$\/-type as indicated,
and of central character in $\Omega(\bG,K^0)$
which are moreover square integrable at
all places where $B$ ramifies. However,
the Ramanujan-Petersson conjecture seems not to be
known in this generality.

\noindent
2.\ Fortunately, when $B$ is ramified at some finite 
place $v$ (which always happens if $[F:\QQ]$ is odd),
then the Ramanujan-Petersson conjecture 
{\em is}\/ known. 
If $F = \QQ$ the automorphic form can be 
lifted to $\GL_2(\QQ)$, where the 
Ramanujan-Petersson conjecture was proved 
by Deligne to follow from the Weil conjectures
(\cite{Del}), which he subsequently proved
(\cite{Del1}). Else, let $B'/F$ be a quaternion 
algebra ramified precisely at the places
where $B$ is except for $\infty_1$ and $v$.
Then our form can be lifted to $B'$, and 
Carayol \cite{Car} showed that they occur in 
the cohomology of a local system on a
Shimura curve with good reduction at $v$.
The Ramanujan-Petersson conjecture then
follows from the Weil-Deligne bounds \cite{Del2}.

\noindent
3.\ More cases of the Ramanujan-Petersson
conjecture can undoubtedly be proved by 
the techniques of \cite{BL} or of \cite{BR}.
The scope of results that this might yield 
is unclear.
}
\end{remarks}
As a result, we get the following
\begin{theorem}
\label{main}
In the notation of Theorem~{\rm\ref{Ram}} assume in
addition that $B$ is ramified in at least one finite
place, and that all the $k_i$\/'s are of the same
parity. Then $\cL^{\vk,-\vk/2}$ is Ramanujan on $X$.
\end{theorem}

\section{Explicit arithmetic examples}
Let  $p_j^{f_j}$ be $g$ powers
of rational primes $p_j$, not necessarily distinct,
and put $r_j = p_j^{f_j} + 1$. Our results suffice 
to give the following
\begin{theorem}
\label{examples}
There exist infinitely many irreducible Ramanujan local 
systems over infinitely many
irreducible $(r_1,\dots,r_g)$\/-regular complexes.
\end{theorem}
\begin{pf}
Let $\Tr_j$ be an
$r_j$\/-regular tree and put $\Tr = \prod_j \Tr_{j}$.
Let $d$ be an integer $\geq$
\[ \max_{\set{p \text{ a rational prime}}}
    \sum_{\set{j|\,p_j = p}} f_j\,.\]
There exists a totally real number field $F$
of degree $d$ over $\QQ$, and pairwise distinct
finite places $v_j$ of $F$ whose residue fields
have $p_j^{f_j}$ elements. There exists a
totally definite quaternion algebra $B/F$ which
is unramified at all the $v_j$\/'s and ramified
over at least one finite prime (again, this is
automatic if $d$ is odd). Viewing each
$\Tr_j$ as the tree associated to
$B_{v_j}^\times \simeq \GL_2(F_{v_j})$,
we get an action of $B^\times$ modulo
its center $F^\times$ on $\Tr$.

Choose an order $\cM$
in $B$, which for simplicity we take to contain
$\cO_F$, and set $S = \set{v_1,\dots,v_g}$. Let
$\cO_{F,S}$ and $\cM_S$ be the localizations
at $S$, namely the elements of $F$ and $B$
that are integral outside of $S$. For an ideal
$N$ of $\cO_F$ prime to $S$ let $\Gamma(N)$ be
the principal congruence subgroup in $\cM_S$,
namely the kernel of the (surjective) reduction
$\cM_S^\times \ra (\cM/N\cM)^\times$ modulo $N$.
Suppose that $N$ is sufficiently small. Then
$\Gamma(N)$ divided by its center acts freely
on $\Tr$. Now take $\vk = (k_1,\dots,k_g)$, where
the $k_j$\/'s are nonnegative integers of the
same parity. Then
$\cL = \Gamma(N) \backslash
                       (\Tr \times V^{\vk,-\vk/2})$
is a local system on $X(N) = \Gamma(N) \backslash \Tr$
if the condition in Theorem \ref{Ram}(2) is satisfied.
In our case it means that if the parity of the
$k_i$\/'s is odd, then there are no elements in
$\cO_{F,S}^\times $ congruent to $1$ modulo $N$
whose norm to $\QQ$ is negative. For example, this
holds if either the $k_i$\/'s are even, or if $-1$
is not in the subgroup of $(\ZZ/(N\cap\ZZ))^\times$
generated by the $p_j^{f_j}$\/'s.

Assuming this, $\cL$ is Ramanujan local system
on $X(N)$ by Theorem~\ref{main}. If we fix the
quaternion algebra and vary the level $N$, we
thus get an infinite family of complexes of the
same regularities with growing number of vertices,
proving the theorem.
\end{pf}

In certain cases the construction takes a
particularly simple form, in which, among other
things, the $\Gr_j$\/'s are all Cayley graphs on
the same group. (The general case is not much
worse --- see \cite[Section 2.6]{JL6} for the
prototypical examples.) Namely, let us assume
there is an ideal $N_0\neq 0$ of $\cO_F$,
prime to the $v_j$\/'s (we allow $N_0 = \cO_F$),
such that the following holds:
\begin{conditions}
\label{specase}
\begin{enumerate}
\item Every ideal of $F$ has a totally positive
generator $\equiv 1\mod N_0$.
\item The class number of $B$ is $1$.
\item The units $\cM^\times$ of a maximal
order $\cM$ of $B$ surject onto
$(\cM/N_0\cM)^\times$, with the kernel being
contained in the center $\cO_F^\times$ of
$\cM^\times$.
\end{enumerate}
\end{conditions}
We then have the following two propositions:
\begin{proposition}
\label{generators}
{\rm 1.}\ For each $1\leq j \leq g$ there are exactly $r_j$
{\rm (}principal\/{\rm )} ideals $P_{j,i}$, $1\leq i \leq r_j$, of
$\cM$ whose norm to $F$ is the prime ideal $v_j$. This
ideal has a totally positive generator, say $\pi_j$,
which is $\equiv 1\mod{N_0}$.  We can next choose
generators \mbox{$\vpi_{j,i}\equiv 1 \mod{N_0\cM}$} for
$P_{j,i}$ whose norm is $\pi_j$, whose image
$\vpi_{j,i}^\ast$ under the main involution
$x\mapsto x^\ast$ of $B$ is some $\vpi_{j,i'}$
for some $1\leq i' \leq r_j$
{\rm (}$i' = i$ may happen{\rm )}.

\noindent
{\rm 2.}\ For every permutation $\sigma$ of $\set{1,\dots,g}$
and any sequence of indices $i_1,\dots,i_g$, with
$1\leq i_j \leq r_j$, there is a
{\rm (}unique\/{\rm )} sequence
$i'_1,\dots, i'_g$, with $1\leq i'_j \leq r_j$, and a
{\rm (}unique\/{\rm )} unit $u\in \cO_F^\times$, satisfying
$u\equiv 1\mod{N_0}$, so that
\begin{equation}
\label{permut}
\vpi_{\sig(1),i_1}\cdots \vpi_{\sig(g),i_g} =
         u\vpi_{1,i'_1}\cdots \vpi_{g,i'_g}\,.
\end{equation}
\end{proposition}
\begin{pf}
We shall use the notations of Section~\ref{qua}.
Fix generators $\pi_j$ of $v_j$ as required. Let
$\cM$ be a maximal order in $B$. For each finite
prime $v$ of $F$ let $\cM_v$ be the completion of
$\cM$ at $v$, and set $K_v = \cM_v^\times$ and
$K^0 = \prod K_v$, the product taken over the
finite $v$\/'s not among the $v_j$\/'s.  Fix the
identifications in (\ref{iden}) so that
$\cM_{v_j} = \Matwo(\cO_{F,v_j})$. The assumption
that $F$ and $B$ have class number one implies,
using the Eichler-Kneser strong approximation
theorem,  that the complex $Y(K^0)$ of (\ref{adel})
coincides with the complex
$X = \Gamma(K^0)\backslash \Tr$ and that it
has one vertex. Therefore there are elements
$\vpi'_{j,i} \in \Gamma(K^0)$, for $1\leq j \leq g$
and $1\leq i \leq r_j$, mapping the vertex of $\Tr$
fixed by $K$ to its $r_j$ neighbours of direction
$\{j\}$. Multiplying by an element in the center, we
may assume that these elements are in $\cM$, and
not divisible by any of the $v_j$. Then the norm
of each $\vpi'_{j,i}$ must be a generator of $v_j$.
Multiplying by a unit in $\cM$ we may assume that
$\vpi'_{j,i} \equiv 1 \mod{N_0\cM}$. Multiplying by
a unit in the center $\cO_F$ we may further assume
that, in addition, the norm is $\pi_j$. This is
because there are units in $\cO_F$, which have
arbitrary signs at the infinite places of $F$ and
which are $\equiv 1\mod{N_0}$. The result gives
the required $\vpi_{j,i}$\/'s. For the second part,
we see from the action on the complex that
$\vpi_{\sig(1),i_1}\cdots \vpi_{\sig(g),i_g} =
         \vpi_{1,i'_1}\cdots \vpi_{g,i'_g}u$
for some unit $u$ of $\cM$. Then $u$ must be
$\equiv 1\mod{N_0}$, and therefore in the center.
\end{pf}

To formulate the second proposition we need to set
some notation first. For a subset $J$ of
$\set{1,\dots,g}$ we will denote by
$\prod_{j\in J}g_j$ the product $g_{j_1}\dots g_{j_n}$
where the $g$\/'s are elements in any semigroup and
the $j_i$\/'s are the elements of $J$ in increasing
order.

Let $N_1$ be a prime ideal of $\cO_F$ prime to the
$v_j$\/'s and to $N_0$, and set $N = N_0N_1$. Let $A$
be the subgroup of $(\cO_F/N_1)^\times$ generated
by the images modulo $N_1$ of the $\pi_j$\/'s. Let
$B$ be the subgroup of scalars in
$(\cM/N_1\cM)^\times$ generated by the images modulo
$N_1\cM$ of the $\pi_j$\/'s and by those units of
$\cO_F$ which are congruent to $1$ modulo $N_0$.
Let $B'$ be the subgroup of scalars in
$(\cM/N_1\cM)^\times$ generated the images modulo
$N_1\cM$ of the $\pi_j$\/'s and by those units of
$\cO_F$ which are congruent to $1$ modulo $N_0$ and
whose norm to $\QQ$ is positive (namely $1$).
The groups
\mbox{$H = \set{g\in (\cM/N_1\cM)^\times \mid \Nm(g) \in A}
       / B$} and
\mbox{$H' = \set{g\in (\cM/N_1\cM)^\times \mid \Nm(g) \in A}
       / B'$}
are isomorphic to the groups
$\SL_2(\cO_F/N_1)$, $\PSL_2(\cO_F/N_1)$, or
$\PGL_2(\cO_F/N_1)$. The latter case occurs
if and only if at least one of the $\pi_j$\/'s is a
square modulo $N_1$. Moreover, $H$ is a
quotient of $H'$, and  the kernel subgroup has order
$1$ if $-1$ is in the subgroup of
$(\ZZ/(N_0\cap \ZZ))^\times$ generated by the norms
\mbox{$\Nm_{F/\QQ} \pi_j \mod{N_0 \cap \ZZ}$}, and has order
$2$ otherwise. If the kernel is $1$ assume that the
parity of the $k_i$\/'s is even. For
$x\in\cM_S^\times$ let $\ov{x}$ and $\ov{x}{}'$
denote its respective reductions into $H$ and $H'$.
These make sense because $\cM_S/N_0\cM_S$ is
isomorphic to $\cM/N_0\cM$. Choose a section
$s:H\ra H'$. With this notation we have the
following

\begin{proposition}
\label{explicit}
{\rm 1}.\ The vertices of $X(N)$ are the elements of
$H$.

\noindent
{\rm 2}.\ The {\rm (}oriented\/{\rm )} cubes of direction
$J\subset \set{1,\dots,g}$ of $X(N)$ are the pairs
$(v,I)$, where $v\in H = \Ver X(N)$ and $I$ is a
$|J|$\/-tuple $(i_j)_{j\in J}$, with
$1\leq i_j \leq r_j$.

\noindent
{\rm 3}. For $j \in J\subset{1,\dots,g}$, the
$j$\/th bottom of the $J$\/-cube $(v,I)$ is
the $J'$\/-cube $(v\ov{\vpi_{j,i'_j}},I')$, where
$J'  = J-\set{j}$, and $I' = (i'_k)_{k\in J'}$
is the $|J'| = (|J|-1)$\/-tuple characterized by
\[\prod_{k\in J} \vpi_{k,i_k} = u\vpi_{j,i_j}
\prod_{k\in J'} \vpi_{k,i'_k}\,,\]
with u a unit $u\in \cO_F^\times$ as in
Proposition {\rm \ref{generators}(2)}.

With the same notation, $\inv{j}(v,I)$ is the
$J$\/-cube $(v\ov{\vpi_{j,i'_j}}, I'')$, with $I''$
characterized by
\[\vpi_{j,i_j}^\ast\prod_{k\in I'} \vpi_{k,i'_k} =
  u'\prod_{k\in J} \vpi_{k,i''_k}   \,\]
for an appropriate $u' \in \cO_F^\times$.

\noindent
{\rm 4}.\ The fiber of $\cL$ over each $v\in H$
is $\cL(v) = V^{\vk,-\vk/2}$. To describe the
transition maps, let $e$ be an edge of direction $\{j\}$,
say $e = (v,i) (= (v,\set{i}))$, with
$1\leq i \leq r_j$. Then $\botj e = v$ and
$\topj e = v\ov{\vpi_{j,i}}$. Set $\epsilon = 1$
if $s(\topj e) = s(v) \ov{\vpi_{j,i}}{}'$ and
$\epsilon = -1$ otherwise. Then
$\cL(e): \cL(\botj e) \ra \cL(\topj e)$ is given
by
\begin{equation}
\label{lc}
\epsilon^{k_1}\rho^{\vk,-\vk/2}(\vpi_{j,i})\,.
\end{equation}
\end{proposition}
\begin{pf}
Let $K^0(N)$ be the principal congruence subgroup
of level $N$ is $K^0$. Using the Eichler-Kneser
strong approximation theorem we see that
the complex $Y(K^0(N))$ is described as
$\Gamma(K^0) \backslash (\Tr \times K^0/CK^0(N))$,
where $C$ is the center of $\Gamma(K^0)$. Then
$X(N)$ is the connected component of $Y(K^0(N))$
lying under $\Tr\times \set{1}$. On the other hand,
$K^0/CK^0(N) \simeq (\cM/N\cM)^\times/B = H$; and
since $\Gamma(K^0)$ acts transitively on the set
\mbox{$\Ver \Tr \times (\cM/N_0\cM)^\times$}, part 1)
follows readily. The rest is routine.
We will show directly (without using $\Tr$) that
the transition maps $\cL (e)$ as in 4) satisfy
the flatness condition
(\ref{flat}). First notice that
$v' = \topf{j_2}\topf{j_1}\sigma$ indeed coincides
with $\topf{j_1}\topf{j_2}\sigma$ for any $2$\/-cell
$\sigma = (v,i_{j_1},i_{j_2})$, with
$1\leq i_j \leq r_j$ for $j= j_1,j_2$. Hence
the lifts of $v$ and of $v'$ are related by
\[s(v')
   = \epsilon' \ov{\vpi_{j_1,i_{j_1}}
              \vpi_{j_2,i_{j_2}}}\,{}' s(v)
   = \epsilon'' \ov{\vpi_{j_2,i'_{j_2}}
          \vpi_{j_1,i'_{j_1}}}\,{}'\tilde{v}\,,
\]
where the product of the signs corresponding to
$\epsilon'$  and $\epsilon''$ is the image in $H$
of $u$ from equation~(\ref{permut}).
Condition~(\ref{flat})
is an immediate consequence.
\end{pf}

The {\em girth} of a cubical complex is the length,
namely the number of edges, of the shortest
homotopically nontrivial closed path.
Set $q = \max_j p_j^{f_j}$ and $n = \Nm_{F/\QQ} N_1$.
We now have the following generalization of
\cite[Theorem 7.3.7]{Lub}
\begin{proposition}
Assume that the group $H$ of $X(N)$ is
$\SL_2(\cO_F/N_1)$. Then the girth of
$\cX(N)$ is at least $2\log_q(N^2/4^d)$.
Put differently,
\[\mbox{{\rm girth} } \cX(N) \geq \frac{4}{3}\log_q
  \# \Ver \cX(N) - \text{{\rm constant}}\,.\]
\end{proposition}
\begin{pf}
As in the $1$\/-dimensional case, the girth is
$\min\dist(x,\gamma x)$, the minimum
taken over all vertices $x$ of the universal cover
$\Tr$ and noncentral elements $\gamma$ of
$\Gamma = \Gamma(K^0)$.
Fix such $\gamma$ and $x$ where the minimum occurs.
Then the distance is, by definition, the sum of
the distances $\dist_j$ of the projections to the
$j$\/th tree factor for $1\leq j \leq g$. Moreover,
$\dist_j \geq -2 \val_{v_j} \Tra_{B/F} \gamma$,
exactly as in loc. cit., Lemma 7.3.2. We have
$N_1^2 |\Nm_{B/F} (\gamma - 1) = 2 - \Tra_{B/F}\gamma$.
Also, since $G_{\infty_i} \simeq \SU(2)$
we have $|\Tra_{B/F}\gamma|_{\infty_i} \leq 2$
for all $1\leq i \leq d$. Taking norms to $\QQ$
gives that $\Nm_{F/\QQ} (2 - \Tra_{B/F}\gamma)$ is
of the form $n^2 m /\prod_j p_j^{l_jf_j}$ for some
$m$, necessarily $\geq 1$, where
$l_j = -\val_{v_j} (2 - \Tra_{B/F}\gamma)$. Hence
\[|\Nm_{F/\QQ} (2 - \Tra_{B/F}\gamma)|
 \leq \prod_{1\leq i \leq d}
        (2 + |\Tra_{B/F}\gamma|_{\infty_i}) \leq 4^d\,,\]
so that
$n^2/4 \leq \prod_j p_j^{l_jf_j} \leq
  q^{\sum_j l_j}.$
Therefore,
\begin{equation}
\begin{array}{rcl}
\text{{\rm girth} } & = & \sum_j \dist_j
  \geq -2\sum_j \min(0,\val_{v_j}\Tra_{F/\QQ}\gamma)
     = 2\sum_{l_j>0} l_j \geq 2\sum_j l_j\\
 & \geq & 2\log_q\left(n^2 / 4^d\right)
 \end{array}
\end{equation}

On the other hand, the cardinality of $\Ver \cX(N)$
is at least that of $\PSL_2(\cO_F/N_1)$, which is at
least a positive constant $c$ times $n^3$. This
implies the last assertion.

\end{pf}
\begin{examples}  {\rm
In \cite[Chapter 5]{Vig} one finds examples satisfying
 Conditions \ref{specase}. These suffice to give
many types of examples. However the regularities are
generally assumed to satisfy some congruence conditions
besides of being primes or prime powers of a specific
type. To get rid of such restrictions it appears necessary
to consider the general cases, whose finite description is
somewhat messier and hence not discussed here.

\noindent
{\bf A.}\ Take $F = \QQ$ and let $B_2$ have discriminant
$2$. When $S$ consists of $g = 1$ prime 
\mbox{$p\equiv 1\mod{4}$},
one gets the Lubotzky-Phillips-Sarnak graphs \cite{LPS}.
The nontrivial local systems over them are handled in 
\cite{JL6}. Notice that the sign $\epsilon$ in
(\ref{lc}) is mistakenly missing there. The 
case of $g = 2$ distinct
primes $p\equiv q \equiv 1 \mod{4}$, without
local systems, is studied for different reasons
in \cite{Moz}.  The resulting complexes are
$(p+1,q+1)$\/-regular.)

\noindent
{\bf B.}\ The case $F = \QQ$, with $B_3$ of discriminant
$3$ and $N_0 = 2$, satisfies  the Conditions
\ref{specase}.

\noindent
{\bf C.}\ The case $F = \QQ$, with $B_{13}$ of discriminant
$13$ and $N_0 = 1$, satisfies  the Conditions
\ref{specase}.
It is used to get $3$\/-regular Ramanujan graphs in
\cite{Chi}.

\noindent
{\bf D.}\ To get regularities other than $p+1$ 
one must use fields other than $\QQ$. The 
simplest case occurs when $F = \QQ(\sqrt{5}\,)$ 
and $B_F = B_2 \otimes F$, with $B_2$ as in A, 
and $N_1 = 2\cO_F$. Then $B_F$ is totally definite
but unramified at all the finite places. If we 
assume the Ramanujan-Petersson conjecture for the 
relevant automorphic forms as in Theorem~\ref{Ram}, 
then this gives irreducible Ramanujan local 
systems over infinitely many $p^2 +1$\/-regular 
graphs for all primes $p\equiv\pm 2 \mod{5}$. 
It also gives irreducible Ramanujan local systems 
over infinitely many $(p+1,p+1)$\/-regular 
complexes for all primes $p\equiv\pm 1 \mod{5}$. 
Of course one can combine regularities to get 
higher dimensional examples. Unfortunately,
Theorem~\ref{main} does not apply and we have
not found an adequate reference in the literature
to prove the Ramanujan-Petersson conjecture 
in this case.

\noindent
{\bf E.}\ For our last example take $F = \QQ(\cos 2\pi/7)$,
and set $B_F = B_2 \otimes F$, with $B_2$ as before.
Then $B_F$ is ramified precisely at the prime $2$ 
of $F$ and at the $3$ infinite places of $F$. Here
Theorem~\ref{main} applies. This example allows
therefore to get infinitely many Ramanujan graphs 
of regularity $p^3+1$ (and Ramanujan local systems 
on them) for any $p$ which is not a square modulo $7$; 
see \cite{Mor} for more general examples coming from
function fields (with the trivial local system). 
Our present example also gives infinitely many 
$(p+1,p+1)$\/-regular and $(p+1,p+1,p+1)$\/-regular
Ramanujan complexes and Ramanujan local systems on
them for any $p$ which is a square modulo $7$. Again,
different such regularities may be combined to form
higher dimensional examples.   }
\end{examples}

\end{document}